# NOTE SUR LE CHOIX DES COURBES FAIT PAR AL-KHAYYAM DANS SA RESOLUTION DES EQUATIONS CUBIQUES ET COMPARAISON AVEC LA METHODE DE DESCARTES

**Nicolas Farès**
Département des mathématiques Faculté des Sciences - Université Libanaise
*Équipe d'Étude et de Recherche sur la Tradition Scientifique Arabe* – CNRS - Liban
nfares55@hotmail.com

***RESUME***

*Dans le Traité algébrique d'al-Khayyâm, se trouve formulée pour la première fois une théorie de résolution des équations cubiques par l'intersection de courbes géométriques. De plus, ce mathématicien a réussi à résoudre, par cette méthode, les quatorze types de ces équations. Il a résolu chacun de ces types au moyen d'un couple de sections coniques prises parmi les suivantes : cercle, parabole et hyperbole (équilatère). Il ne justifie pas son choix du couple de courbes utilisées dans la solution : son style en ce qui concerne ce point est purement synthétique.*

*Cette note porte sur un détail pointu : justifier le choix fait par al-Khayyâm de ces couples de courbes, bien que ce détail ait été étudié par R. Rashed. Par la même occasion, nous nous sommes posé la même question pour Descartes (dont le style n'est pas moins synthétique), à la recherche d'une ressemblance possible des deux méthodes de résolution.*

*Nous avons pu constater que le choix des quatorze couples de courbes a été fait par al-Khayyâm conformément à un procédé de calcul uniforme qu'il appliquait systématiquement à chacun des quatorze types d'équation.*

*Le choix de tels couples de courbes chez Descartes semble être le résultat de techniques de calcul bien distinctes, bien que le projet soit, sur le fond, le même et que la motivation soit, dans l'essentiel, la même: résolution par intersection de coniques utilisant leurs équations.*

**Mots-clés:** Histoire de la géométrie algébrique. Résolutions géométriques des équations cubiques. Al-Khayyâm - Descartes.

***ABSTRACT***

*It is well known that Al-Khayâm, for the first time in history, formulated a complete theory to solve third degree equations using the intersection of geometric curves and moreover solved the fourteen types of equations using this method. His solution for these equations was either using the intersection of two parabolas, a parabola and a circle, a parabola and a hyperbola, a circle and a hyperbola or two hyperbolas. His style was purely synthetic lacking any analysis which may lead to deducing any motives for his choice of the curves.*

*Our article is centered on a pointed detail which is the justification of the choice of these curves, knowing that R. Rashed studied this choice in all of the mentioned fourteen equations. As for Descartes, whose style is as synthetic as Al-Khayâm's, we ask the same question concerning the solution of the third degree equations, searching for any possible resemblance between the two methods of solution, both mathematicians having been motivated similarly in geometric resolutions.*

*We noticed that the choice of the fourteen curves couples was made by Al Khayyâm regarding a uniform calculation technique that he applied systematically to each equation type among the fourteen.*

*It has appeared that the choice of the curves made by Descartes is the result of distinct calculation techniques, even though the projects of the two mathematicians are of the same essence, and their motivations are the same: finding the solution of third degree equations using the intersection of conical curves and the equations of curves.*
**Keywords:** History of the geometrical algebra. Geometrical resolutions of the cubic equations. Al-Khayyâm – Descartes.

**1. Introduction**

La rédaction de cette note est motivée par la traduction en arabe du livre intitulé *al-Khayyâm mathématicien,* œuvre récente de R. Rashed, avec la participation de B. Vahabzadeh (R. Rashed et B. Vahabzadeh , 1999) . Ce livre contient tous les travaux mathématiques d'Omar al-Khayyâm parvenus à nos jours, édités, traduits et présentés avec des paraphrases et des commentaires mathématiques. Ces travaux d'al-Khayyâm peuvent être subdivisés en deux parties.

La première et la plus volumineuse occupe la première partie du livre ; dans cette partie, se trouve formulée, pour la première fois dans l'histoire, une théorie de résolution des équations cubiques par l'intersection de courbes géométriques ; al-Khayyâm l'expose et résout ensuite effectivement les types d'équations algébriques de degré inférieur ou égal à 3. Cette partie se situe dans le courant *géométrique*[1] du développement de l'algèbre.

La deuxième partie de ces travaux s'attaque d'une façon critique aux fondements des mathématiques ; elle contient trois traités:
- un premier Traité qui a pour but de "démontrer" le cinquième postulat d'Euclide, projet qui a occupé les mathématiciens sans interruption pendant plus de vingt siècles et qui a donné naissance aux géométries non euclidiennes,
- deux traités sur les rapports et la proportionnalité, soulignant les difficultés des définitions euclidiennes des livres V et VI des *Eléments* d'Euclide. Les textes des traités de cette deuxième partie doivent intéresser particulièrement les philosophes et les chercheurs en philosophie des mathématiques.

La publication de ce livre constitue, nous semble-t-il, la dernière étape de la réalisation d'un projet de R. Rashed qui est l'écriture de l'histoire de l'algèbre arabe. Il a en effet exposé l'histoire de l'algèbre « *arithmétique* » dans ses deux livres : *Al- Bâhir fi'l-jabr d'al-Samaw'al* (S. Ahmad et R. Rashed, 1972 ) et *Entre arithmétique et algèbre, Recherche sur l'histoire des mathématiques arabes* (R. Rashed, 1984). Il a commencé la rédaction de l'histoire de l'autre courant, le courant géométrique, du développement de l'algèbre, dans son livre *Sharaf al-Dîn al-Tûsî : œuvres mathématiques - Algèbre et géométrie au 12ème siècle* (R. Rashed, 1986). Dans ce livre, R. Rashed promet, entre autres choses, de faire une étude comparative de l'algèbre géométrique d'al-Khayyâm et de celle de Descartes. Il fallait bien attendre 13 ans pour que cette promesse soit tenue. L'introduction du livre *al-Khayyâm mathématicien* est, en effet, un article profond et difficile dans lequel R. Rashed présente une étude inédite de l'histoire de l'algèbre *géométrique*. Il aboutit à l'évaluation des éléments de modernité dans la *Géométrie* de Descartes en se basant sur les nouvelles données que présentent les travaux mathématiques de la tradition arabe dont les travaux d'al-Khayyâm (qui occupent la première partie du livre susmentionné) et de Sharaf al-Dîn al-Tûsî constituent l'apogée[2].

---

[2]On sait, en effet, depuis l'œuvre de R. Rashed *Entre arithmétique et algèbre* (R. Rashed, 1984) que cette nouvelle discipline mathématique (l'algèbre) s'est développée depuis sa fondation par al-khwârizmî dans deux courants : *arithmétique* et *géométrique. Le* courant *arithmétique* a été étudié dans l'œuvre susmentionné ainsi que dans une deuxième œuvre de R. Rashed et S. Ahmed : "*al- Bâhir en algèbre d'al-ssamaw'al*" (S. Ahmed et R. Rashed, 1972) . Le courant géométrique est étudié dans l'introduction d'un troisième livre de R. Rashed: "*Sharaf al-Dîn al-Tûsî -Œuvres mathématiques – algèbre et géométrie au XIIe siècle*" (R. Rashed, 1986).

Notre présente note se base exclusivement sur cette première partie du livre ; elle porte sur un détail pointu, à savoir, justifier le choix des courbes fait par al-Khayyâm dans sa résolution des équations cubiques, bien que ce détail ait été étudié par R. Rashed. En effet, pour chaque type d'équation, al-Khayyâm utilise un couple de courbes dont l'intersection donne la solution, mais il ne justifie pas son choix de ces courbes : son style en ce qui concerne ce point est purement synthétique. Nous ne nous attaquons ici qu'au côté relatif aux détails des calculs qui l'ont conduit à ce choix. Par la même occasion, nous nous sommes posé la même question pour Descartes (dont le style n'est pas moins synthétique), à la recherche d'une ressemblance possible. Or, il semble, comme nous allons le voir, que les deux types de choix utilisent des techniques de calcul distinctes, bien que le projet est, sur le fond, le même comme l'a montré l'étude de R. Rashed et que les motivations sont en gros les mêmes : résolution par intersection de coniques.

**2. Choix des courbes fait par al-Khayyâm.**

Al-Khayyâm a divisé les équations cubiques en quatorze types :

Equation 3 : $x^3 = c$. ………………………….. …... pp.
Equation 13 : $x^3 + bx = c$. …………… ……...pc.
Equation 14 : $x^3 + c = bx$. …………… ……...ph.
Equation 15 : $x^3 = bx + c$. …………… ……...ph.
Equation 16 : $x^3 + ax^2 = c$. …………… ……...ph.
Equation 17 : $x^3 + c = ax^2$. …………… ……...ph.
Equation 18 : $x^3 = ax^2 + c$. ………………. .ph.
Equation 19 : $x^3 + ax^2 + bx = c$. ……………ch.
Equation 20 : $x^3 + ax^2 + c = bx$. ……………hh.
Equation 21 : $x^3 + bx + c = ax^2$. ……………ch.
Equation 22 : $x^3 = ax^2 + bx + c$. ……………hh.
Equation 23 : $x^3 + ax^2 = bx + c$. ……………hh.
Equation 24 : $x^3 + bx = ax^2 + c$. ……………ch.
Equation 25 : $x^3 + c = ax^2 + bx$. ……………hh.

Il résout chacun de ces type par l'intersection de deux coniques prises parmi les suivantes: le cercle, la parabole et l'hyperbole (équilatère) [3].

R. Rashed a interprété le choix des courbes pour chacune des équations et, pour la compréhension de ce choix des courbes il a conseillé de reprendre le raisonnement d'al-Khayyâm en sens inverse (R. Rashed et B. Vahabzadeh , 1999, p. 37) : partant de la fin vers le début. L'application minutieuse de ce conseil nous a permis de remarquer qu'il y a une seule idée derrière le choix de ces couples de courbes. Cette idée principale réside dans le lemme permettant de construire la racine carrée d'une quantité algébrique, ce qui revient à construire la moyenne proportionnelle entre deux quantités algébriques.

Les grandes lignes de la démarche d'al-Khayyâm peuvent se résumer pratiquement par le schéma suivant : l'équation est transformée en l'égalité de deux rapports :

---

[2] Dans cette tradition commencée avec al-Khwârizmî, on peut citer de noms illustres: Thâbit Ibn Qurra, al-Mâhâni, al- Bîrûni, Abû al-Jûd Ibn al-Layth, Abû Nasr Ibn Irâq, Abû Jafar al- Khâzin, …,

[3] La numérotation de ces équations est celle adoptée par R. Rashed conformément à l'ordre de leur traitement dans le Traité d'al-Khayyâm. (Les deux lettres qui suivent l'expression d'une équation indiquent les courbes qui ont été utilisées dans sa résolution : *pp* = deux paraboles, *pc* = parabole et cercle, *h* = hyperbole équilatère)

$\dfrac{M_1}{M_2} = \dfrac{N_1}{N_2}$, où $M_1$ et $M_2$ sont deux surfaces (de dimension 2) et $N_1$ et $N_2$ sont de dimension 1. La moyenne proportionnelle $y$ de $N_1$ et de $N_2$ donne alors la racine carrée de $\dfrac{M_1}{M_2}$ et, de la relation : $\dfrac{\sqrt{M_1}}{\sqrt{M_2}} = \dfrac{N_1}{y} = \dfrac{y}{N_2}$, il tire les équations des courbes qui résolvent l'équation en question.

Pour illustrer cette démarche nous allons décrire (dans un langage moderne) la "marche en sens inverse" pour chacune de ces équations. Par souci de concision, nous croyons qu'il suffit de décrire la démarche d'al-Khayyâm seulement pour les équations 3, 13 et 14. Nous remarquons que les équations où $b \neq 0$ sont traitées toutes, en utilisant le rapport $\dfrac{b}{x^2} = \left(\sqrt{\dfrac{b}{x}}\right)^2$ que, plus haut, nous avons noté $\dfrac{M_1}{M_2}$. Nous allons commencer par ces équations (en nombre de huit) laissant à la fin les 4 autres (équations 16, 17, 18 et 3), où $b = 0$.

**Equation 13 : $x^3 + bx = c$**

Démarche d'al-Khayyâm : Cette démarche peut se résumer comme suit :

Soient les segments de droite AB = $\sqrt{b}$, $h$ = BC, perpendiculaire à AB tel que $b.h = c$, (i.e. $h = \dfrac{c}{b}$). Considérons la parabole (P) de sommet B, d'axe AB, de côté droit AB (i.e. $y\sqrt{b} = x^2$) et le cercle (C) de diamètre BC, (i.e. $y^2 = x.(\dfrac{c}{b} - x)$). Soit D leur point d'intersection et soient G, la projection de D sur AB et E sa projection sur BC et soit $x$ = DG = EB, (l'abscisse de D) et $y$ = DE = BG (son ordonnée).

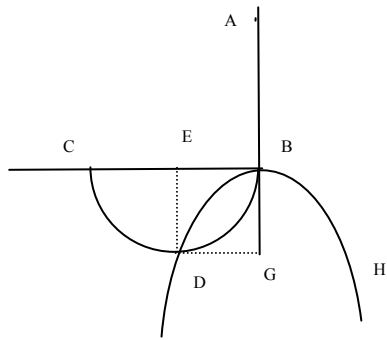

On a : $y\sqrt{b} = x^2$, car D∈(P), d'où $\dfrac{\sqrt{b}}{x} = \dfrac{x}{y}$, mais D appartient aussi à (C) d'où $\dfrac{x}{y} = \dfrac{y}{\frac{c}{b} - x}$, d'où $\dfrac{\sqrt{b}}{x} = \dfrac{x}{y} = \dfrac{y}{\frac{c}{b} - x}$, d'où $\dfrac{b}{x^2} = \dfrac{x}{\frac{c}{b} - x}$ ($= \dfrac{x}{y} \times \dfrac{y}{\frac{c}{b} - x}$), d'où $b \cdot (\dfrac{c}{b} - x) = x^3$, d'où $x^3 = c - bx$, d'où ...$x^3 + bx = c$ et, $x$ est racine de l'équation.

<u>Marche en sens inverse</u> : L'équation 13 : $x^3 + bx = c$ s'écrit $x^3 = c - bx$, c. à. d. $b \cdot (\dfrac{c}{b} - x) = x^3 = x^2 \cdot x$, d'où $\dfrac{b}{x^2} = \dfrac{x}{\frac{c}{b} - x}$, i.e. $\left(\dfrac{\sqrt{b}}{x}\right)^2 = \dfrac{x}{\frac{c}{b} - x}$. Si on prend la moyenne proportionnelle $y$ entre $x$ et $\dfrac{c}{b} - x$, c. à. d. tel que : $\dfrac{x}{y} = \dfrac{y}{\frac{c}{b} - x}$, on a

(E) $\dfrac{\sqrt{b}}{x} = \dfrac{x}{y} = \dfrac{y}{\frac{c}{b} - x}$,

d'où le choix des courbes fait par al-Khayyâm. En effet, l'équation (E) peut être remplacée par deux quelconques des 3 équations suivantes :

(1) $\dfrac{\sqrt{b}}{x} = \dfrac{x}{y}$, (qui donne l'équation de (P)),

(2) $\dfrac{x}{y} = \dfrac{y}{\frac{c}{b} - x}$, (qui donne l'équation de (C)),

(3) $\dfrac{\sqrt{b}}{x} = \dfrac{y}{\frac{c}{b} - x}$, (hyperbole équilatère).

**Remarque 1 :** Al-Khayyâm a choisi le couple des courbes donné par les équations (1) et (2). Il avait deux autres choix convenables offerts par les couples d'équations: ((1), (3)) et ((2), (3)). Nous ne savons pas pourquoi il a préféré le choix qu'il a fait parmi les trois offerts par cette méthode. De même, dans chacune des autres équations (14 – 25), al-Khayyâm n'utilise qu'un seul choix, sans signaler l'existence des deux autres choix offerts par cette méthode.

**Equation 14** : $x^3 + c = bx$

Démarche d'al-Khayyâm : Cette démarche peut se résumer comme suit :

Soient les segments de droite AB = $\sqrt{b}$, $h$ = BC, perpendiculaire à AB tel que $b.h = c$ (i.e. $h = \dfrac{c}{b}$). Considérons la parabole (P) de sommet B, d'axe AB et de côté droit AB (i.e. $y\sqrt{b} = x^2$) et l'hyperbole équilatère (H) de sommet C, d'axe BC et de côté droit BC (i.e. $y^2 = x.(x - \dfrac{c}{b})$).

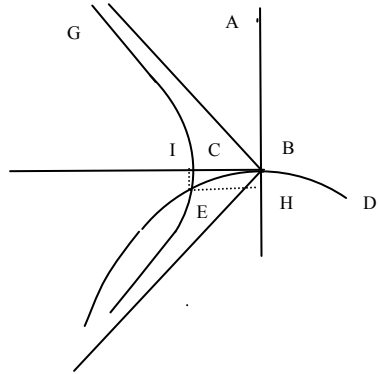

Al-Khayyâm remarque que si ces deux courbes ne se coupent pas le problème est impossible ; mais si elle se coupent, il considère un de leurs deux points d'intersection, soit E. Soient I, la projection de E sur BC et H sa projection sur AB et soient $x$ = BI = EH (l'abscisse de E) et, $y$ = EI = BH (son ordonnée).

On a : E∈(H), d'où $y^2 = x.(x - \dfrac{c}{b})$, d'où $\dfrac{x}{y} = \dfrac{y}{x - \dfrac{c}{b}}$. Mais, E∈(P), donc $y\sqrt{b} = x^2$

d'où $\dfrac{\sqrt{b}}{x} = \dfrac{x}{y}$, d'où $\dfrac{\sqrt{b}}{x} = \dfrac{x}{y} = \dfrac{y}{x - \dfrac{c}{b}}$, d'où $\dfrac{b}{x^2} = \dfrac{x}{x - \dfrac{c}{b}}$ (= $\dfrac{x}{y} \times \dfrac{y}{x - \dfrac{c}{b}}$), d'où

$b \cdot (x - \dfrac{c}{b}) = x^3$, d'où $x^3 + c = bx$ et, $x$ est racine de l'équation.

<u>Marche en sens inverse</u> : L'équation 14 : $x^3 + c = bx$ s'écrit $x^3 = bx - c$, c. à. d. $b \cdot (x - \frac{c}{b}) = x^3 = x^2 \cdot x$, d'où $\frac{b}{x^2} = \frac{x}{x - \frac{c}{b}}$, i.e. $\left(\frac{\sqrt{b}}{x}\right)^2 = \frac{x}{x - \frac{c}{b}}$. Si on prend la moyenne proportionnelle $y$ entre $x$ et $x - \frac{c}{b}$, c. à. d. tel que : $\frac{x}{y} = \frac{y}{x - \frac{c}{b}}$, on a

(E) $\quad \dfrac{\sqrt{b}}{x} = \dfrac{x}{y} = \dfrac{y}{x - \frac{c}{b}}$,

d'où le choix des courbes fait par al- Khayyâm. En effet, l'Équation (E) peut être remplacée par deux quelconques des 3 équations suivantes :

(1) $\dfrac{\sqrt{b}}{x} = \dfrac{x}{y}$, (qui donne l'équation de (P)),

(2) $\dfrac{x}{y} = \dfrac{y}{x - \frac{c}{b}}$, (qui donne l'équation de (H)).

(3) $\dfrac{\sqrt{b}}{x} = \dfrac{y}{x - \frac{c}{b}}$, (hyperbole équilatère).

Al-Khayyâm a choisi le couple des courbes, donné par les équations (1) et (2). (Voir la remarque1 plus haut).

Comme nous l'avons souligné, à partir de l'équation suivante, nous allons décrire avec concision, seulement *la marche en sens inverse* du raisonnement d'al-Khayyâm.

**Equation 15** : $x^3 = bx + c$

$$x^3 = bx + c \Leftrightarrow x^2 \cdot x = b \cdot (x + \frac{c}{b}) \Leftrightarrow \frac{b}{x^2} = \frac{x}{x + \frac{c}{b}} \Leftrightarrow \left(\frac{\sqrt{b}}{x}\right)^2 = \frac{x}{x + \frac{c}{b}}.$$

Donc, si $y$ est la moyenne proportionnelle entre $x$ et $x + \frac{c}{b}$, c. à. d. tel que : $\frac{x}{y} = \frac{y}{x + \frac{c}{b}}$,

on a

$$(E) \quad \frac{\sqrt{b}}{x} = \frac{x}{y} = \frac{y}{x+\dfrac{c}{b}},$$

d'où le choix des courbes fait par al-Khayyâm. En effet, l'Équation (E) peut être remplacée par deux quelconques des 3 équations suivantes :

(1) $\dfrac{\sqrt{b}}{x} = \dfrac{x}{y}$, (la parabole (P)),

(2) $\dfrac{x}{y} = \dfrac{y}{x+\dfrac{c}{b}}$, (l' hyperbole (H)),

(3) $\dfrac{\sqrt{b}}{x} = \dfrac{y}{x+\dfrac{c}{b}}$, (hyperbole équilatère).

Al-Khayyâm a choisi le couple des courbes donné par les équations (1) et (2). (Voir la remarque1, plus haut).

**Equation 19 :** $x^3 + ax^2 + bx = c.$

$$x^3 + ax^2 + bx = c \Leftrightarrow x^3 + ax^2 = c - bx \Leftrightarrow x^2.(x+a) = b.\left(\frac{c}{b} - x\right)$$

$$\Leftrightarrow \frac{b}{x^2} = \frac{x+a}{\dfrac{c}{b} - x} \Leftrightarrow \left(\frac{\sqrt{b}}{x}\right)^2 = \frac{x+a}{\dfrac{c}{b} - x}.$$

Donc, si $y$ est la moyenne proportionnelle entre $x+a$ et $\dfrac{c}{b} - x$, c. à. d. tel que :

$\dfrac{x+a}{y} = \dfrac{y}{\dfrac{c}{b} - x}$, on a

$$(E) \quad \frac{\sqrt{b}}{x} = \frac{x+a}{y} = \frac{y}{\dfrac{c}{b} - x},$$

d'où le choix des courbes fait par al-Khayyâm. En effet, l'Équation (E) peut être remplacée par deux quelconques des 3 équations suivantes :

(1) $\dfrac{\sqrt{b}}{x} = \dfrac{y}{\dfrac{c}{b} - x}$, (l'hyperbole (H)),

(2) $\dfrac{x+a}{y} = \dfrac{y}{\dfrac{c}{b} - x}$, (le cercle (C)),

(3) $\dfrac{\sqrt{b}}{x} = \dfrac{x+a}{y}$, (parabole).

Al-Khayyâm a choisi le couple des courbes donné par les équations (1) et (2). (Voir la remarque1 plus haut).

**Equation 20** : $x^3 + ax^2 + c = bx$

$$x^3 + ax^2 + c = bx \Leftrightarrow x^3 + ax^2 = bx - c \Leftrightarrow x^2.(x+a) = b.\left(x - \dfrac{c}{b}\right)$$

$$\Leftrightarrow \dfrac{b}{x^2} = \dfrac{x+a}{x - \dfrac{c}{b}} \Leftrightarrow \left(\dfrac{\sqrt{b}}{x}\right)^2 = \dfrac{x+a}{x - \dfrac{c}{b}}.$$

Donc, si $y$ est la moyenne proportionnelle entre $x+a$ et $x - \dfrac{c}{b}$, c. à. d. tel que :

$\dfrac{x+a}{y} = \dfrac{y}{x - \dfrac{c}{b}}$, on a

$$(E)\quad \dfrac{\sqrt{b}}{x} = \dfrac{x+a}{y} = \dfrac{y}{x - \dfrac{c}{b}},$$

d'où le choix des courbes fait par al-Khayyâm En effet, l'Équation (E) peut être remplacée par deux quelconques des 3 équations suivantes :

(1) $\dfrac{\sqrt{b}}{x} = \dfrac{y}{x - \dfrac{c}{b}}$, (l'hyperbole ($H_1$)),

(2) $\dfrac{x+a}{y} = \dfrac{y}{x - \dfrac{c}{b}}$, (l'hyperbole ($H_2$)),

(3) $\dfrac{\sqrt{b}}{x} = \dfrac{x+a}{y}$, (parabole).

Al-Khayyâm a choisi le couple des courbes donné par les équations (1) et (2). (Voir la remarque1 plus haut).

**Equation 21 :** $x^3 + bx + c = ax^2$.

$$x^3 + bx + c = ax^2 \Leftrightarrow ax^2 - x^3 = bx + c \Leftrightarrow x^2.(a-x) = b.(x+\dfrac{c}{b})$$

$$\Leftrightarrow \dfrac{b}{x^2} = \dfrac{a-x}{x+\dfrac{c}{b}} \Leftrightarrow \left(\dfrac{\sqrt{b}}{x}\right)^2 = \dfrac{a-x}{x+\dfrac{c}{b}}.$$

Donc, si $y$ est la moyenne proportionnelle entre $a - x$ et $x + \dfrac{c}{b}$, c. à. d. tel que :

$\dfrac{a-x}{y} = \dfrac{y}{x+\dfrac{c}{b}}$, on a

$$(E)\quad \dfrac{\sqrt{b}}{x} = \dfrac{a-x}{y} = \dfrac{y}{x+\dfrac{c}{b}},$$

d'où le choix des courbes fait par al- Khayyâm En effet, l'Équation (E) peut être remplacée par deux quelconques des 3 équations suivantes :

(1) $\dfrac{\sqrt{b}}{x} = \dfrac{y}{x+\dfrac{c}{b}}$, (l'hyperbole (H)),

(2) $\dfrac{a-x}{y} = \dfrac{y}{x+\dfrac{c}{b}}$, (le cercle (C)),

(3) $\dfrac{\sqrt{b}}{x} = \dfrac{a-x}{y}$, ( parabole).

Al-Khayyâm a choisi le couple des courbes donné par les équations (1) et (2). (Voir la remarque1 plus haut).

**Equation 22** : $x^3 = ax^2 + bx + c.$

$$x^3 = ax^2 + bx + c \Leftrightarrow x^3 - ax^2 = bx + c \Leftrightarrow x^2.(x-a) = b.(x+\frac{c}{b})$$

$$\Leftrightarrow \frac{b}{x^2} = \frac{x-a}{x+\frac{c}{b}} \Leftrightarrow \left(\frac{\sqrt{b}}{x}\right)^2 = \frac{x-a}{x+\frac{c}{b}}.$$

Donc, si $y$ est la moyenne proportionnelle entre $x - a$ et $(x+\frac{c}{b})$, c. à. d. tel que :

$\dfrac{x-a}{y} = \dfrac{y}{x+\dfrac{c}{b}}$ , on a

$$(E) \quad \frac{\sqrt{b}}{x} = \frac{x-a}{y} = \frac{y}{x+\frac{c}{b}},$$

d'où le choix des courbes fait par al- Khayyâm. En effet, l'Équation (E) peut être remplacée par deux quelconques des 3 équations suivantes :

(1) $\dfrac{\sqrt{b}}{x} = \dfrac{y}{x+\dfrac{c}{b}}$ , (l'hyperbole ($H_1$)),

(2) $\dfrac{x-a}{y} = \dfrac{y}{x+\dfrac{c}{b}}$ , (l'hyperbole ($H_2$)),

(3) $\dfrac{\sqrt{b}}{x} = \dfrac{x-a}{y}$ , (parabole).

Al-Khayyâm a choisi le couple des courbes donné par les équations (1) et (2). (Voir la remarque1 plus haut).

**Equation 23** : $x^3 + ax^2 = bx + c.$

$$x^3 + ax^2 = bx + \Leftrightarrow x^2.(x+a) = b.(x+\frac{c}{b})$$

$$\Leftrightarrow \frac{b}{x^2} = \frac{x+a}{x+\frac{c}{b}} \Leftrightarrow \left(\frac{\sqrt{b}}{x}\right)^2 = \frac{x+a}{x+\frac{c}{b}}.$$

Donc, si $y$ est la moyenne proportionnelle entre $x + a$ et $(x + \frac{c}{b})$, c. à. d. tel que :

$\frac{x+a}{y} = \frac{y}{x+\frac{c}{b}}$, on a

$$(E) \quad \frac{\sqrt{b}}{x} = \frac{x+a}{y} = \frac{y}{x+\frac{c}{b}},$$

d'où le choix des courbes fait par al-Khayyâm. En effet, l'Équation (E) peut être remplacée par deux quelconques des 3 équations suivantes :

(1) $\frac{\sqrt{b}}{x} = \frac{y}{x+\frac{c}{b}}$, (l'hyperbole ($H_1$)),

(2) $\frac{x+a}{y} = \frac{y}{x+\frac{c}{b}}$, (l'hyperbole ($H_2$)),

(3) $\frac{\sqrt{b}}{x} = \frac{x+a}{y}$, (parabole).

Al-Khayyâm a choisi le couple des courbes donné par les équations (1) et (2). (Voir la remarque1, plus haut).

**Equation 24 :** $x^3 + bx = ax^2 + c$

$$x^3 + bx = ax^2 + c \Leftrightarrow ax^2 - x^3 = bx - c \Leftrightarrow x^2.(a-x) = b.(x - \frac{c}{b})$$

$$\Leftrightarrow \frac{b}{x^2} = \frac{a-x}{x-\frac{c}{b}} \Leftrightarrow \left(\frac{\sqrt{b}}{x}\right)^2 = \frac{a-x}{x-\frac{c}{b}}.$$

Donc, si $y$ est la moyenne proportionnelle entre $a - x$ et $(x - \frac{c}{b})$, c. à. d. tel que :

$\frac{a-x}{y} = \frac{y}{x-\frac{c}{b}}$, on a

$$(E) \quad \frac{\sqrt{b}}{x} = \frac{a-x}{y} = \frac{y}{x - \frac{c}{b}},$$

d'où le choix des courbes fait par al- Khayyâm. En effet, l'Équation (E) peut être remplacée par deux quelconques des 3 équations suivantes :

(1) $\dfrac{\sqrt{b}}{x} = \dfrac{y}{x - \dfrac{c}{b}}$, (l'hyperbole (H)),

(2) $\dfrac{a-x}{y} = \dfrac{y}{x - \dfrac{c}{b}}$, (le cercle (C)),

(3) $\dfrac{\sqrt{b}}{x} = \dfrac{a-x}{y}$, (parabole).

Al-Khayyâm a choisi le couple des courbes donné par les équations (1) et (2). (Voir la remarque1 plus haut).

**Equation 25** : $x^3 + c = ax^2 + bx$.

$$x^3 + c = ax^2 + bx \Leftrightarrow x^3 - ax^2 = bx - c \Leftrightarrow x^2.(x-a) = b.(x - \frac{c}{b})$$

$$\Leftrightarrow \frac{b}{x^2} = \frac{x-a}{x - \frac{c}{b}} \Leftrightarrow \left(\frac{\sqrt{b}}{x}\right)^2 = \frac{x-a}{x - \frac{c}{b}}.$$

Donc, si $y$ est la moyenne proportionnelle entre $x - a$ et $(x - \dfrac{c}{b})$, c. à. d. tel que :

$\dfrac{x-a}{y} = \dfrac{y}{x - \dfrac{c}{b}}$, on a

$$(E) \quad \frac{\sqrt{b}}{x} = \frac{x-a}{y} = \frac{y}{x - \frac{c}{b}},$$

d'où le choix des courbes fait par al- Khayyâm. En effet, l'Équation (E) peut être remplacée par deux quelconques des 3 équations suivantes :

(1) $\dfrac{\sqrt{b}}{x} = \dfrac{y}{x - \dfrac{c}{b}}$, (l'hyperbole ($H_1$)),

(2) $\dfrac{x - a}{y} = \dfrac{y}{x - \dfrac{c}{b}}$, (l'hyperbole ($H_2$)),

(3) $\dfrac{\sqrt{b}}{x} = \dfrac{x - a}{y}$, (parabole).

Al-Khayyâm a choisi le couple des courbes donné par les équations (1) et (2). (Voir la remarque1 plus haut).

**Equation 16** : $x^3 + ax^2 = c$.

$$x^3 + ax^2 = c \Leftrightarrow x^2 \cdot (x + a) = c^{\frac{1}{3}} \cdot c^{\frac{2}{3}}$$

$$\Leftrightarrow \dfrac{c^{\frac{2}{3}}}{x^2} = \dfrac{x + a}{c^{\frac{1}{3}}} \quad \left( = \left( \dfrac{c^{\frac{1}{3}}}{x} \right)^2 \right)$$

Donc, si $y$ est la moyenne proportionnelle entre $x + a$ et $c^{\frac{1}{3}}$, c. à. d. tel que : $\dfrac{x + a}{y} = \dfrac{y}{c^{\frac{1}{3}}}$,

on a

$$(E) \quad \dfrac{c^{\frac{1}{3}}}{x} = \dfrac{x + a}{y} = \dfrac{y}{c^{\frac{1}{3}}},$$

d'où le choix des courbes fait par al-Khayyâm. En effet, l'Équation (E) peut être remplacée par deux quelconques des 3 équations suivantes :

(1) $\dfrac{c^{\frac{1}{3}}}{x} = \dfrac{y}{c^{\frac{1}{3}}}$, (l'hyperbole (H)),

(2) $\dfrac{x + a}{y} = \dfrac{y}{c^{\frac{1}{3}}}$, (la parabole (P)),

(3) $\dfrac{c^{\frac{1}{3}}}{x} = \dfrac{x+a}{y}$, (parabole).

Al-Khayyâm a choisi le couple des courbes donné par les équations (1) et (2). (Voir la remarque1 plus haut).

**Equation 17 :** $x^3 + c = ax^2$.

$$x^3 + c = ax^2 \Leftrightarrow x^2.(a-x) = c^{\frac{1}{3}} \cdot c^{\frac{2}{3}}$$

$$\Leftrightarrow \dfrac{c^{\frac{2}{3}}}{x^2} = \dfrac{a-x}{c^{\frac{1}{3}}} \quad \left(= \left(\dfrac{c^{\frac{1}{3}}}{x}\right)^2\right)$$

Donc, si $y$ est la moyenne proportionnelle entre $a - x$ et $c^{\frac{1}{3}}$, c. à. d. tel que : $\dfrac{a-x}{y} = \dfrac{y}{c^{\frac{1}{3}}}$,

on a

(E) $\dfrac{c^{\frac{1}{3}}}{x} = \dfrac{a-x}{y} = \dfrac{y}{c^{\frac{1}{3}}}$

d'où le choix des courbes fait par al- Khayyâm. En effet, l'Équation (E) peut être remplacée par deux quelconques des 3 équations suivantes :

(1) $\dfrac{c^{\frac{1}{3}}}{x} = \dfrac{y}{c^{\frac{1}{3}}}$, (l'hyperbole (H)),

(2) $\dfrac{a-x}{y} = \dfrac{y}{c^{\frac{1}{3}}}$, (la parabole (P)),

(3) $\dfrac{c^{\frac{1}{3}}}{x} = \dfrac{a-x}{y}$, (parabole).

Al-Khayyâm a choisi le couple des courbes donné par les équations (1) et (2). (Voir la remarque1 plus haut).

**Equation 18 :** $x^3 = ax^2 + c$.

$$x^3 = ax^2 + c \Leftrightarrow x^2.(x-a) = c = \dfrac{c}{a} \cdot a$$

$$\Leftrightarrow \frac{\left(\dfrac{c}{a}\right)}{x^2} = \frac{x-a}{a} \Leftrightarrow \left(= \left(\frac{\sqrt{\dfrac{c}{a}}}{x}\right)^2\right)$$

Donc, si *y* est la moyenne proportionnelle entre *x - a* et *a*, c. à. d. tel que :

$$\frac{x-a}{y} = \frac{y}{a}$$

on a :

$$(E) \quad \frac{\sqrt{\dfrac{c}{a}}}{x} = \frac{x-a}{y} = \frac{y}{a},$$

d'où le choix des courbes fait par al- Khayyâm. En effet, l'Équation (E) peut être remplacée par deux quelconques des 3 équations suivantes :

(1) $\dfrac{\sqrt{\dfrac{c}{a}}}{x} = \dfrac{y}{a}$, (l'hyperbole (H)),

(2) $\dfrac{x-a}{y} = \dfrac{y}{a}$, (la parabole (P)),

(3) $\dfrac{\sqrt{\dfrac{c}{a}}}{x} = \dfrac{x-a}{y}$, (parabole).

(Al-Khayyâm a choisi le couple des courbes donné par les équations (1) et (2). (Voir la remarque1 plus haut).

**Remarque2:** Nous ne savons pas pourquoi al-Khayyâm. n'a pas utilisé la même technique des équations 16 et 17, en écrivant $c = c^{\frac{1}{3}} \cdot c^{\frac{2}{3}}$, *ce qui aurait conduit à*

$$x^3 = ax^2 + c \Leftrightarrow x^2 \cdot (x-a) = c^{\frac{1}{3}} \cdot c^{\frac{2}{3}}$$

$$\Leftrightarrow \frac{c^{\frac{2}{3}}}{x^2} = \frac{x-a}{c^{\frac{1}{3}}} \quad \left(= \left(\frac{c^{\frac{1}{3}}}{x}\right)^2\right)$$

Alors, si *y* est la moyenne proportionnelle entre *x - a* et $c^{\frac{1}{3}}$, c. à. d. tel que :

$$\frac{x-a}{y} = \frac{y}{c^{\frac{1}{3}}}$$

on a

$$(E) \quad \frac{c^{\frac{1}{3}}}{x} = \frac{x-a}{y} = \frac{y}{c^{\frac{1}{3}}},$$

d'où un choix de couple de courbes semblable aux choix correspondant dans les équations 16 et 17.

D'une autre part, al-Khayyâm aurait pu utiliser pour les équations 16 et 17 la même méthode de l'équation 18 utilisant la décomposition : $c = \frac{c}{a} \cdot a$.

**Equation 3 :** $x^3 = c$.

Cette équation a un statut spécial. Al-Khayyâm ne la résout pas directement par l'intersection de coniques mais, à l'aide d'un lemme déjà démontré, il prend entre *1* et *c*, [4] deux moyennes proportionnelles, *x* et *y* (telles que $\frac{1}{x} = \frac{x}{y} = \frac{y}{c}$), alors *x* est la racine cherchée. Or, c'est la démonstration de ce lemme (R. Rashed et B. Vahabzadeh , 1999 , p. 31, lemme1]) : *construction de deux segments x et y entre deux segments a et c tels que les quatre soient en proportion continue*, i.e tels que:

$$(E) \quad \frac{a}{x} = \frac{x}{y} = \frac{y}{c}$$

qui a été faite par l'intersection de deux paraboles : $x^2 = a.y$ et $y^2 = c.x$.

Comme l'a remarqué R. Rashed *le choix des courbes s'explique par la proportion notée* c. à. d. par les deux équations constituant la relation *(E)* :

$$\frac{a}{x} = \frac{x}{y}, (x^2 = a.y) \quad \text{et} \quad \frac{x}{y} = \frac{y}{c}, (y^2 = c.x)$$

On peut bien remarquer qu'une troisième conique aurait pu être choisie pour la solution (pour remplacer une de ces paraboles), qui est l'hyperbole équilatère *x.y = a.c*, donnée par l'équation $\frac{a}{x} = \frac{y}{c}$.

Si l'on suit le raisonnement d'al-Khayyâm en sens inverse, on remarque que l'équation 3 s'écrit : $x^2.x = 1^2 .c,$ d'où

$$\frac{1^2}{x^2} = \frac{x}{c}$$

---

[4] Il s'agit du segment unité et d'un segment de longueur *c* (hauteur du solide égal à *c* et dont la base est un carré de côté 1.

alors $\dfrac{1}{x} = \sqrt{\dfrac{x}{c}}$ et, si l'on prend la moyenne proportionnelle y entre *x* et *c*, on a

$\dfrac{x}{y} = \dfrac{y}{c} = \sqrt{\dfrac{x}{c}}$, d'où $\dfrac{1}{x} = \dfrac{x}{y} = \dfrac{y}{c}$

ce qui pourrait justifier le recours au *lemme 1* susmentionné.

**3. Choix des courbes fait par Descartes.**
3.1. Description sommaire de la méthode de Descartes.
  Descartes résout simultanément l'équation du 3$^e$ degré et l'équation du 4$^e$ degré à l'aide d'un cercle et d'une parabole (Ch. Adam et P. Tannery, 1982, pp. 464-468).
  Il rappelle que toute équation de degré n
$$x^n + a_{n-1}.x^{n-1} + a_{n-2}.x^{n-2} + \ldots a_1.x + a_0$$
s'écrit sous la forme
$$x^n + b_{n-2}.x^{n-2} + b_{n-3}.x^{n-3} + \ldots b_1.x + b_0$$

(il réduit à zéro le coefficient $a_{n-1}$ de $x^{n-1}$, par un changement de variable $z = x + \dfrac{a}{n}$ (Ch. Adam et P. Tannery, 1982, p. 449) ; alors l'équation générale du 3$^e$ degré prend la forme :
(1)        $z^3 \propto *.\ apz\,.\ aaq$ .
Il la résout en considérant l'équation du 4$^e$ degré
(2)        $z^4 \propto *.\ apzz\,.\ aaqz\,.\ a^3r$ ;
où :
  $\propto$ désigne le signe de l'égalité ;
  * désigne la place vide d'une puissance de z ;
  **.** (le point) désigne l'un des signes + ou – (les divers cas correspondant aux valeurs, + ou – du point « **.** », correspondent aux divers types d'équations) **;**
  *a* désigne l'unité, écrite pour le respect formel de l'homogénéité.
  Notons que, dans la suite de son travail, il se passe de l'écriture de a et les équations (1) et (2) sont écrites respectivement sous la forme :
(1')        $z^3 \propto *.\ pz\,.\ q$
(2')        $z^4 \propto *.\ pzz\,.\ qz\,.\ r$
que nous pouvons écrire, avec des notations modernes :
{1}        $z^3 + pz + q = 0$,
{2}        $z^4 + pz^2 + qz + r = 0$.
  Sa méthode de résolution simultanée est purement synthétique. L'exposé de cette méthode n'en facilite pas la compréhension à cause de la multitude des cas entremêlés d'une façon encombrante[5]. Il est pourtant clair que Descartes considère (sans le dire explicitement) l'équation {1} comme cas particulier de {2} où $r = 0$, $(z^4 + pz^2 + qz = 0 \Leftrightarrow z^3 + pz + q = 0)$ ; il suffit donc de résoudre {2}.
  Pour cela il prend la parabole (P) : $y = z^2$, qu'il utilise pour tous les cas de figure. Les solutions de {2} sont alors données par l'intersection de cette parabole et de divers cercles choisis selon les cas de figure. En exposant ces différents cas Descartes utilise des segments de longueurs $\frac{1}{2}, \frac{1}{2}p, \frac{1}{2}q, \frac{1}{2}p - \frac{1}{2}$, sans aucune explication.

---
[5] Suivant que l'un ou l'autre des coefficients *p, q, r*, soit >, < ou = 0.

La présence de ces segments est difficilement justifiable sans la reconstitution de l'analyse qui a conduit au raisonnement de Descartes. En nous basant notamment sur le calcul qu'il donne à la fin de son exposé pour vérifier, dans un des cas de figure (Ch. Adam et P. Tannery, 1982, pp. 467-68), que le segment désignés est bien racine de l'équation {2}, il nous semble que cette analyse peut être reconstituée de la façon suivante :

En posant $y = z^2$, on a :

$$z^4 + pz^2 + qz + r = 0 \Leftrightarrow (z^4 + pz^2 - z^2) + (z^2 + qz) + r = 0$$

$$\Leftrightarrow \left(z^2 + \frac{(p-1)}{2}\right)^2 + \left(z + \frac{q}{2}\right)^2 + r - \frac{(p-1)^2}{4} - \frac{q^2}{4} = 0$$

$$\Leftrightarrow \left(y + \frac{(p-1)}{2}\right)^2 + \left(z + \frac{q}{2}\right)^2 - R^2 = 0,$$

cette dernière équation étant bien celle de la famille des cercles utilisées par Descartes. Nous remarquons que dans le cas d'une équation du troisième degré, $R^2$ est $> 0$ car $r = 0$, et l'équation a toujours une solution (réelle).

**4. Quelques remarques sur les techniques utilisées.**

Nous allons nous limiter ici, à souligner quelques détails qui pourraient aider à comparer les méthodes de résolution des équations du 3$^e$ degré chez al-Khayyâm et Descartes.
1) Les deux méthodes visent explicitement la résolution de ces équations géométriquement, par l'intersection de coniques. Al-Khayyâm adopte cette méthode en remarquant que ni lui ni ses prédécesseurs n'ont pas pu le faire par radicaux (R. Rashed et B. Vahabzadeh, 1999 p. 125) tout en exprimant le souhait que « *d'autres, qui nous succèderont sauront-ils le faire* ». Descartes l'adopte délibérément, en soulignant l'insuffisance des méthodes par radicaux, faites par les italiens (Ch. Adam et P. Tannery, 1982, pp. 472 et 474). Cette option pour la méthode géométrique chez Descartes fait partie de son projet plus large visant la résolution des équations algébriques de, tout degré, par intersection de courbes. Les deux méthodes se basent sur les équation des coniques. Mais, alors qu'al-Khayyâm utilise les propriétés des coniques qui reviennent à leurs équations, et que la notion d'équation de ces courbes n'était chez lui qu'implicite, la théorie des équations de courbes de degrés quelconques, était fondamentale dans *la Géométrie* de Descartes.

Un détail pourrait être significatif sur ce plan : Descartes résout, même les équations du second degré, conformément à cette méthode (intersection d'une droite et d'un cercle), tandis qu'al-Khayyâm résout ces équations par la méthode classique d'al-Khwârizmî, (sans citer ce dernier, mais en citant le livre II d'Euclide), méthode qui repose sur l'égalité des surfaces rectangulaires correspondant à chaque type d'équation. Cela pourrait s'expliquer par l'intention d'al-Khayyâm de focaliser son attention sur le problème difficile qui n'a pas été formulé ni résolu par ces prédécesseurs, à savoir la résolution des équations du troisième degré.
2) Un deuxième détail est significatif : dans l'équation cubique
$x^3 + ax^2 + bx + c = 0,$

l'inconnue $x$ et ses puissances $x^2$ et $x^3$ sont chez al-Khayyâm, représentables respectivement par des segments de droite, des carrés et des cubes ; pour le respect de l'homogénéité, $a$ est représentable par un segment de droite et $b$ par une aire (rectangulaire). Or, toutes ces grandeurs sont pour Descartes des « lignes toutes simples », c.à.d. des segments de droites (Ch. Adam et P. Tannery, 1982, p. 370-371). Cela marque un pas important vers la négligence de l'homogénéité donc vers l'indépendance de l'algèbre et de son développement. Notons toutefois que, chez les algébristes arabe, depuis al-Khwârizmî, en passant par ses successeurs, que ça soit du « courant arithmétique » ou du « courant géométrique », l'homogénéité est explicitement dépassée, surtout dans la formulation des équations algébriques. Notons aussi que Descartes même, n'a pas pu s'en débarrasser complètement : « *comme s'il faut tirer la racine cubique de aabb –b, il faut penser que la quantité aabb est divisée une fois par l'unité & que l'autre quantité b est multipliée deux fois par la même* » (Ch. Adam et P. Tannery, 1982, pp. 371-372).

Cette ligne unité a été utilisée dans le calcul géométrique d'al-Khayyâm, Elle a joué, avec son carré et son cube, un rôle important dans la théorie de résolution géométrique des équations algébriques de degré ≤ 3 fondée par lui. On la rencontre chez ses prédécesseurs, Banû Mûsa (IX$^e$ s.) et Ibn al-Haytham (XI$^e$ s.), comme le dit R. Rashed (R. Rashed et B. Vahabzadeh, 1999, p. 8). Plus tard, le carré de cette unité et son cube ont pris chez Sharaf al-Dîn al-Tûsî, les noms d'*unité plane* et d'*unité solide*, respectivement.

3) Il est à noter que le symbolisme algébrique moderne (la notation de l'inconnue et de ses puissances, $x, x^2, .., ax, bx^2, ..$, les signes des opérations : +, -, /, $\sqrt[2]{}$, $\sqrt[3]{}$ -noté $\sqrt{C...}$ -), absent dans la tradition arabe, mais déjà mûr avec Descartes, a beaucoup aidé la réalisation de son projet. Il en est de même pour l'utilisation des racines négatives (racines *fausses ou moindres que rien*) et des coefficients négatifs (dotés du signe -) bien que cette utilisation n'était qu'à son début[6].

4) Descartes a, de plus, utilisé une technique efficace, celle de changement de variable affine ($x \to x - \dfrac{a}{n}$) pour réduire la forme d'une équation, chose qui n'a pas été pratiquée par al-Khayyâm et qui aurait réduit le nombre de types d'équations cubiques qu'il a traités. Cette même technique sera pourtant habilement utilisée par al-Tûsî, successeur d'al-Khayyâm. Notons toutefois, qu'al-Khayyâm lui-même utilise des changements de variables du type $x \to \dfrac{1}{x}$.

Enfin, soulignons un détail qui ne manque pas d'attirer l'attention du lecteur qui feuillette la Géométrie de Descartes et le Traité algébrique d'al-Khayyâm : curieusement, les axes *cartésiens* apparaissent sur les figures utilisées par al-Khayyâm, mais, seul l'axe des ordonnées apparaît explicitement sur celles de Descartes.

**5. Conclusion.**

Rappelons que la comparaison des deux projets mathématiques, celui d'al-Khayyâm et celui de Descartes et leurs apports au fondement de la géométrie algébrique, dépasse le cadre de cet exposé. Ce sujet a été traité par R. Rashed dans l'introduction de son livre susmentionné, en utilisant une documentation et, surtout des connaissances, qui ne sont pas à notre portée. Une œuvre récente de C. Houzel, souligne l'importance des contributions des deux mathématiciens au fondement de la géométrie algébrique (C. Houzel, 2002; voir notamment p.p. 5-14). Nous nous limitons ici, à la minuscule tâche de comparer leurs techniques de résolution des équations cubiques.

---

[6] (Par exemple, $a$ et +$a$, désignaient un coefficient positif, il est désigné par –$a$ quand il est négatif. La lettre $a$, sans le signe –, ne peut pas représenter un coefficient négatif. Cela a constitué un inconvénient, en obligeant Descartes à considérer plusieurs cas dans sa résolution de l'équation cubique qu'il a notée $z^3 \propto *.pz.q$
suivant que les points qui y figurent désignent des signes + ou -. Cela revient bien à considérer plusieurs types d'équations).

D'après notre reconstitution de l'analyse du raisonnement d'al-Khayyâm (§2. plus haut), il s'avère que son choix des courbes de résolution résulte de l'utilisation d'un moyen unique, qu'il applique à tous les types d'équations, à savoir la considération d'une moyenne proportionnelle entre deux segments. La multitude des types d'équation revient à celle des cas correspondant aux signes des coefficients de l'équation.

Du livre de R. Rashed qui se réfère à Beeckmann, nous comprenons qu'en 1619, bien avant la parution de sa *Géométrie* (1637), Descartes a entamé la résolution « *de quatre équations cubiques parmi les treize équations obtenues en combinant nombre, racine, carré et cube.... Sa classification est celle d'al-Khayyâm, ses résultats n'atteignent pas encore la généralité qu'on trouve chez ce dernier,…* » (R. Rashed et B. Vahabzadeh, 1999, pp. 14-15). La dernière formulation de sa méthode de résolution des équations cubiques, qu'il présente dans sa *Géométrie*, utilise une parabole unique et un cercle variable selon le type d'équation. Nous venons de voir (§3. plus haut) que ce choix des courbes repose sur le passage à une équation du 4$^e$ degré qui permet, en posant $y = x^2$, d'obtenir par un calcul simple, l'équation de la famille de cercles en question.

Les deux méthodes reposent donc sur des techniques de calcul distinctes, bien qu'elles sont dictées par un souci et une conception qui sont pratiquement les mêmes et qui d'ailleurs sont liés à un projet qui est historiquement le même : traduire les problèmes de la géométrie en équations algébriques, puis résoudre celles-ci par intersection de courbes géométriques. Si ce projet a été assez mûr et bien formulé par al-Khayyâm, il a ses racines, bien avant lui, dans les travaux de ses prédécesseurs : al- Mahânî, abû-al-Jûd, al-Khâzin ; Ibn-al-Haytham,… . Depuis la naissance de l'algèbre au IX$^e$ siècle, ce projet s'est formé et s'est nourri à partir de problèmes hérités de la mathématique grecque, notamment de la géométrie, et s'est développé en employant les acquis de cette géométrie, celle d'Euclide et surtout celle des coniques d'Apollonius . Avec Descartes il a connu un essor où il a été repris, conçu et formulé dans toute sa généralité. Les équations cubiques, qui y ont joué un rôle important ont été résolues avec des moyens plus développés dont l'usage des quantités négatives mais, surtout, d'un symbolisme algébrique parfait, qui a accéléré le développement de l'algèbre et de toutes les mathématiques. Jusqu'à nos jours, ce projet, objet de la géométrie algébrique, ne cesse de s'animer, usant du développement de toutes les mathématiques et contribuant, à son tour, à leur développement.

**Bibliographie**

<div style="text-align:center">

\*\*\*

\*

</div>